\documentclass[10pt]{amsart}

\usepackage{amsmath,amsthm,amsfonts,latexsym,amssymb,mathrsfs,setspace,enumerate,color,cite,graphicx,epsf,mathtools,fullpage,array,amstext,verbatim,enumitem,float}

\newtheorem{thm}{Theorem}
\newtheorem{lem}[thm]{Lemma}

\numberwithin{equation}{section}
\numberwithin{thm}{section}

\newtheorem{conj}[thm]{Conjecture}

\newcommand{\rat}{\mathbb Q}
\newcommand{\real}{\mathbb R}
\newcommand{\com}{\mathbb C}
\newcommand{\alg}{\overline\rat}
\newcommand{\algt}{\alg^{\times}}

\newcommand{\intg}{\mathbb Z}
\newcommand{\nat}{\mathbb N}

\newcommand{\rank}{\mathrm{rank}}

\newcommand{\xx}{{\bf x}}
\newcommand{\yy}{{\bf y}}
\newcommand{\zz}{{\bf z}}

\newcommand{\ttt}{{\bf t}}

\newcommand{\bb}{{\bf b}}

\newcommand{\tors}{\mathrm{tors}}

\newcolumntype{P}[1]{>{\centering\arraybackslash}p{#1}}

\title[A Polynomial Time Test for Exceptional Points]{A Polynomial Time Test to Detect Numbers with Many Exceptional Points}
\author{Ryan Carpenter \and Charles L. Samuels}

\address{Christopher Newport University, Department of Mathematics, 1 Avenue of the Arts, Newport News, VA 23606}
\email{ryan.carpenter@cnu.edu, charles.samuels@cnu.edu}

\begin{document}

\begin{abstract}
	For each algebraic number $\alpha$ and each positive real number $t$, the $t$-metric Mahler measure $m_t(\alpha)$ creates an extremal problem whose solution varies depending on the value of $t$.  
	The second author studied the points $t$ at which the solution changes, called {\it exceptional points for $\alpha$}.  Although each algebraic number has only finitely many exceptional points, it is conjectured that, 
	for every $N \in \nat$, there exists a
	number having at least $N$ exceptional points.  In this article, we describe a polynomial time algorithm for establishing the existence of numbers with at least $N$ exceptional points.  Our work constitutes an 
	improvement over the best known existing algorithm which requires exponential time.  We apply our main result to show that there exist numbers with at least $37$ exceptional points, another improvement 
	over previous work which was only able to reach $11$ exceptional points.
\end{abstract}

\maketitle

\section{Introduction} \label{Intro}

Suppose $\alpha$ is a non-zero algebraic number of degree $d$ over $\rat$ with minimal polynomial $F$ over $\intg$.  Further suppose that $a$ is the leading coefficient of $F$ and
$\alpha_1,\alpha_2,\ldots,\alpha_d\in \com$ are such that
\begin{equation*}
	F(z) = a\cdot \prod_{i=1}^d (z - \alpha_i).
\end{equation*}
Under these assumptions, the {\it (logarithmic) Mahler measure} of $\alpha$ is defined to be 
\begin{equation*}
	m(\alpha) = \log |a| + \sum_{i=1}^d \log\max\{1,|\alpha_i|\}.
\end{equation*}
It is obvious from the definition that $m(\alpha) \geq 0$ for all $\alpha\in \algt$, and moreover, it follows from Kronecker's Theorem \cite{Kronecker} that $m(\alpha) = 0$ if and only if $\alpha$ is a root of unity.
We also note that the behavior of $m(\alpha)$ is particularly straightforward when $\alpha\in \rat^\times$.  Indeed, if $\alpha = r/s$ and $\gcd(r,s) = 1$ then $m(\alpha) = \log\max\{|r|,|s|\}$.

In attempting to construct large prime numbers, D.H. Lehmer \cite{Lehmer} came across the problem of determining whether there exists a sequence of non-zero
algebraic numbers $\{\alpha_n\}$, not roots of unity, such that $m(\alpha_n)$ tends to $0$ as $n\to\infty$.  This problem remains unresolved, although substantial evidence
suggests that no such sequence exists (see \cite{AmorosoDvornicich, BDM, Schinzel, Smyth}, for instance).  This assertion is usually called Lehmer's conjecture.

\begin{conj}[Lehmer's Conjecture]
	There exists $c>0$ such that $m(\alpha) \geq c$ whenever $\alpha\in \algt$ is not a root of unity. 
\end{conj}

Dobrowolski \cite{Dobrowolski} provided the best known lower bound on $m(\alpha)$ in terms of $\deg\alpha$, while Voutier \cite{Voutier}
later gave a  version of this result with an effective constant.  Nevertheless, little progress has been made on Lehmer's conjecture for an arbitrary algebraic number $\alpha$.

Dubickas and Smyth \cite{DubSmyth, DubSmyth2} were the first to study a modified version of the Mahler measure which gives rise to a metric on $\algt/\algt_\tors$.
Specifically, let
\begin{equation*}
	\mathcal P(\alpha) = \left\{ (\alpha_1,\alpha_2,\ldots,\alpha_n)\in (\algt)^n: n\in \nat,\ \alpha = \prod_{i=1}^n \alpha_i\right\}
\end{equation*}
and define the {\it metric Mahler measure} by
\begin{equation} \label{m1}
	m_1(\alpha) = \inf\left\{ \sum_{i=1}^n m(\alpha_i): (\alpha_1,\alpha_2,\ldots,\alpha_n)\in \mathcal P(\alpha)\right\}.
\end{equation}
It is verified in \cite{DubSmyth2} that $m_1:\algt\to [0,\infty)$ satisfies the following key properties:
\begin{enumerate}[label={(\roman*)}]
	\item\label{Mod} $m_1(\alpha) = m_1(\zeta\alpha)$ for all $\alpha\in \algt$ and $\zeta\in \algt_\tors$
	\item\label{Inverse} $m_1(\alpha) = m_1(\alpha^{-1})$ for all $\alpha \in \algt$
	\item\label{Tri} $m_1(\alpha\beta) \leq m_1(\alpha) + m_1(\beta)$ for all $\alpha,\beta\in \algt$.
\end{enumerate}
These facts combine to ensure that $(\alpha,\beta) \mapsto m_1(\alpha\beta^{-1})$ is a well-defined metric on $\algt/\algt_\tors$ which
induces the discrete topology if and only if Lehmer's conjecture is true.

The second author \cite{SamuelsCollection, SamuelsParametrized,SamuelsMetric} extended the metric Mahler measure to form a parametrized family of metric Mahler measures.
For each real number $t>0$, the {\it $t$-metric Mahler measure} of $\alpha$ is defined to be
\begin{equation*} \label{mt}
	m_t(\alpha) = \inf\left\{ \left(\sum_{i=1}^n m(\alpha_i)^t\right)^{1/t}:  (\alpha_1,\alpha_2,\ldots,\alpha_n)\in \mathcal P(\alpha)\right\}
\end{equation*}
and we note that $m_1(\alpha)$ agrees with the definition provided by Dubickas and Smyth in \eqref{m1}.  Properties \ref{Mod} and \ref{Inverse} continue to hold with $m_t$ in place of $m_1$,
and instead of \ref{Tri}, we have that
\begin{equation*}
	m_t(\alpha\beta)^t \leq m_t(\alpha)^t + m_t(\beta)^t
\end{equation*} 
for all $\alpha,\beta\in \algt$ and all $t > 0$.  As a result, $(\alpha,\beta) \mapsto m_t(\alpha\beta^{-1})^t$ defines a metric on $\algt/\algt_\tors$ which induces the discrete topology if and only 
if Lehmer's conjecture is true.

Although it is known that the infimum in the definition of $m_t(\alpha)$ is attained for all $\alpha\in \algt$ (see \cite{SamuelsInfimum,SamuelsCollection}), the point which attains the infimum depends on $t$.
A positive real number $t$ is called an {\it exceptional point for $\alpha$} if, for every neighborhood $U$ of $t$ and every $(\alpha_1,\alpha_2,\ldots,\alpha_n)\in \mathcal P(\alpha)$, there exists $s\in U$ such that
\begin{equation*}
	m_s(\alpha) < \left( \sum_{i=1}^n m(\alpha_i)^s\right)^{1/s}.
\end{equation*}
Roughly speaking, exceptional points for $\alpha$ are positive real numbers $t$ at which the infimum attaining point in $m_t(\alpha)$ must change.
The second author \cite{SamuelsMetric} established the following fact about the set of exceptional points for a given algebraic number.

\begin{thm} \label{FiniteExceptional}
	Every algebraic number has finitely many exceptional points.
\end{thm}

While Theorem \ref{FiniteExceptional} is a strong result, its proof provides no assistance in listing the exceptional points for a particular algebraic number $\alpha$.  As a result, we obtain an important open question.
Does there exist a uniform upper bound, independent of $\alpha$, on the number of exceptional points for $\alpha$?  Based on the second author's work \cite{SamuelsContinued,SamuelsExceptional},
we suspect that no such bound exists even for rational $\alpha$.

\begin{conj} \label{ExceptionalCounting}
	For every $N\in \nat$ there exists a rational number having at least $N$ exceptional points.
\end{conj}

The second author \cite{SamuelsExceptional} verified Conjecture \ref{ExceptionalCounting} for $N=11$, although unfortunately, our techniques required a brute force computational approach which was unsuccessful 
for $N\geq 12$.  Specifically, for each $n\geq 3$ we defined a certain finite collection $\mathcal V_n$ of $n$-dimensional vectors over $\intg$.  In order to create a number having at least $n-2$ exceptional points, 
we needed to compute and store the complete list of elements in $\mathcal V_n$, and unfortunately, $\#\mathcal V_n$ grows at least exponentially in terms of $n$ (see Theorem \ref{ExpGrowth} for a formal proof of this assertion).  
For these reasons, the techniques of \cite{SamuelsExceptional} should be expected to fail even for relatively small values of $n$.  We were unable to perform the required computations for any $n\geq 14$.

In the present article, we describe a polynomial time algorithm to check whether there exists a rational number having at least $n-2$ exceptional points.
Specifically, let $f_n$ denote the Fibonacci sequence defined so that $f_0=0$, $f_1=1$ and $f_{n} = f_{n-1} + f_{n-2}$ for all $n\geq 2$.
For each pair of integers $(k,\ell)$ with $1\leq k < \ell \leq n$, we define the vector $\yy_n(k,\ell) \in \real^n$ to be given by $\yy_n(k,\ell) = (y_1,y_2,\ldots,y_n)^T$, where
\begin{equation*}
	y_i = \begin{cases}
		(-1)^{\ell - k + 1}f_{n-\ell}/f_{\ell - k} & \mbox{if } i = k\\
		f_{n-k}/f_{\ell - k} & \mbox{if } i = \ell\\
		0 & \mbox{if } i\not\in \{k,\ell\}.
		\end{cases}
\end{equation*}
Although the definition of $\yy_n(k,\ell)$ is technical, it is simple to compute.  For instance, we have the following:
\begin{equation*}
	\yy_7(3,4) = \begin{pmatrix} 0 \\ 0 \\ f_{3}/f_{1} \\ f_{4}/f_{1} \\ 0 \\ 0 \\ 0 \end{pmatrix} = \begin{pmatrix} 0 \\ 0 \\ 2 \\ 3 \\ 0 \\ 0 \\ 0 \end{pmatrix},\
		\yy_9(2,5) = \begin{pmatrix} 0 \\ f_4/f_3 \\ 0 \\ 0 \\ f_7/f_3 \\ 0 \\ 0 \\ 0\\ 0\end{pmatrix} = \begin{pmatrix} 0 \\ 3/2 \\ 0 \\ 0 \\ 13/2 \\ 0 \\ 0 \\0 \\ 0 \end{pmatrix},\ 
		\yy_{11}(4,10) = \begin{pmatrix} 0 \\0 \\ 0 \\ -f_1/f_6 \\0 \\ 0 \\ 0 \\ 0\\ 0 \\ f_7/f_6 \\ 0 \end{pmatrix} =  \begin{pmatrix} 0 \\0 \\ 0 \\ -1/8 \\0 \\ 0 \\ 0 \\ 0\\ 0 \\ 13/8 \\ 0 \end{pmatrix}.
\end{equation*}
Now define
\begin{equation*}
	\mathcal F_n =\{\yy_n(\ell,k):1\leq k < \ell\leq n\}\quad\mbox{and}\quad \mathcal F^+_n = \left\{ \yy_n(k,\ell): 1\leq k < \ell \leq n,\ \ell - k\ \mathrm{is\ odd}\right\}.
\end{equation*}
We note that all entries of $\yy_n(k,\ell)$ are non-negative if and only if either $\ell - k$ is odd or $n=\ell$.  However, in the latter case one easily verifies that $\yy_n(k,\ell) = \yy_n(\ell-1,\ell) \in \mathcal F^+_n$.
Therefore, $\mathcal F_n^+$ is indeed the set of all vectors in $\mathcal F_n$ for which all entries are non-negative.

We let $\phi = (1+\sqrt 5)/2$, and for each $\xx = (x_1,x_2,\ldots,x_n)^T\in \real^n$, we define $G_{\xx}:\real^+ \to \real$ by
\begin{equation*}
	G_{\xx}(t) = \sum_{i=1}^n x_i\max\{f_i,\phi f_{i-1}\}^t.
\end{equation*}
Our main result shows how to use the above information to test, in polynomial time, for the existence of numbers with at least $n-2$ exceptional points.

\begin{thm} \label{Main}
	Let $n\geq 3$ be an integer and $t_1,t_2,\ldots,t_{n-2},t_{n-1}\in \real^+$ be such that $t_1 > t_2 > \cdots > t_{n-2} > t_{n-1}$. Further assume the following:
	\begin{enumerate}[label={(\Alph*)}]
		\item\label{UniquenessMain} For each $1\leq i\leq n-1$ there exists a unique point $\zz(i)\in \mathcal F^+_n$ such that $$G_{\zz(i)}(t_i) = \min\{ G_{\xx}(t_i):\xx\in \mathcal F^+_n\}.$$
		\item\label{IntegersMain} $\zz(i)\in \intg^n$ for all $1\leq i\leq n-1$.
		\item\label{NoMatchMain} $\zz(i) \ne \zz(i+1)$ for all $1\leq i\leq n-2$.
	\end{enumerate}
	Then there exist infinitely many pairs of primes $(p,q)$ such that $p^{f_n}/q^{f_{n-1}}$ has at least $n-2$ exceptional points.
\end{thm}

The assumptions of Theorem \ref{Main} can be verified by calculating $G_{\xx}(t)$ at each element $(\xx,t)\in \mathcal F^+_n\times\mathcal T_n$, where $\mathcal T_n = \{t_1,t_2,\ldots,t_{n-2},t_{n-1}\}$.
An elementary counting argument (see Theorem \ref{FCard}) reveals that
\begin{equation} \label{FTCardinality}
	\#(\mathcal F^+_n\times\mathcal T_n) = \begin{cases} (n^3 -3n^2 + 7n - 5)/4 & \mbox{if } n\mbox{ is odd} \\ (n^3 - 3n^2 + 6n - 4)/4 & \mbox{if } n\mbox{ is even}, \end{cases}
\end{equation}
so as promised, Theorem \ref{Main} creates a polynomial time algorithm to verify the existence of rational numbers having at least $n-2$ exceptional points.
If we wish to apply Theorem \ref{Main} to a particular value of $n$, one crucial challenge remains -- to select $t_1,t_2,\ldots,t_{n-2},t_{n-1}$ in a way that satisfies the required properties.  
Of course, we do not wish to select the $t_i$ arbitrarily and to simply hope for favorable results.  Fortunately, the results of \cite{SamuelsExceptional} help to provide important guidance. 

If $i\geq 1$ is an integer, then according to \cite[Lemma 4.1]{SamuelsExceptional}, there exists a unique positive real number $s_i$ such that
\begin{equation*}
	\max\{f_{i+2},\phi f_{i+1}\}^{s_i} = \max\{f_{i+1},\phi f_{i}\}^{s_i} + \max\{f_{i},\phi f_{i-1}\}^{s_i}.
\end{equation*}
Further observations from \cite[\S 4]{SamuelsExceptional} establish that $\{s_i\}_{i=1}^\infty$ is a strictly decreasing sequence with $s_i>1$ for all $i$.  By using a straightforward continuity argument, it is
possible to show that $s_i$ approaches $1$ as $i\to\infty$. 

Now that we have defined the sequence $\{s_i\}_{i=1}^\infty$, we may select points $t_1,t_2,t_3,\ldots$ such that
\begin{equation*}
	t_1 > s_1 > t_2 > s_2 > t_3 > s_3 > \cdots.
\end{equation*}
Based on evidence presented in \cite{SamuelsExceptional} and Section \ref{Application}, we suspect that $t_1,t_2,\ldots,t_{n-2},t_{n-1}$ always satisfy the assumptions of Theorem \ref{Main}.  
Using choices of $t_i$ of this type, we have verified the assumptions of Theorem \ref{Main} for $n = 39$, leading to the following improvement to the work of \cite{SamuelsExceptional}.

\begin{thm} \label{MainExceptional}
	There exist infinitely many rational numbers having at least $37$ exceptional points.
\end{thm}

The proof of Theorem \ref{MainExceptional} is a computation performed in MATLAB which utilizes the algorithm described above.  Clearly this is a large improvement over \cite{SamuelsExceptional}
which established the existence of numbers with at least $11$ exceptional points.  Unlike the work of \cite{SamuelsExceptional}, our primary computational limitation does not arise from the running time of our algorithm.
Instead, we encounter a machine precision error when $n=40$, in which case MATLAB is unable to distinguish among some values of $G_{\xx}(t_{38})$ even if those values really are distinct 
(see Section \ref{Application} for further detail on why this occurs precisely at $n=40$ but not at $n=39$).  Consequently, we are unable to verify the assumptions of Theorem \ref{Main} when $n=40$, 
and cannot verify the existence of numbers with at least $38$ exceptional points.

It would be an interesting direction of future research to seek improvements to our computational
methods in a way that enables the relevant calculations when $n\geq 40$.  As it currently stands, we find the existing evidence presented here to be compelling support for Conjecture \ref{ExceptionalCounting}.

The remainder of this article is structured in the following way.  We use Section \ref{MainProof} to prove our main result, Theorem \ref{Main}.
In Section \ref{Application}, we use the polynomial search resulting from Theorem \ref{Main} to establish Theorem \ref{MainExceptional}.  
Finally, Section \ref{Supplementary} contains proofs of two supplementary results that are relevant to our work.  

\section{Proof of Theorem \ref{Main}} \label{MainProof}

We shall construct the proof of Theorem \ref{Main} by applying a series of four lemmas.  The first of those lemmas concerns the solution set to a certain matrix equation involving the Fibonacci sequence.

\begin{lem} \label{2x2Lemma}
	If $k, \ell$ and $n$ are integers with $1\leq k< \ell \leq n$ then $x_1 = (-1)^{\ell-k+1}f_{n-\ell}/f_{\ell - k}$ and $x_2 = f_{n-k}/f_{\ell - k}$ defines the unique solution to the equation
	\begin{equation} \label{2x2Solve}
		\begin{pmatrix} f_k & f_\ell \\ f_{k-1} & f_{\ell - 1} \end{pmatrix}\begin{pmatrix} x_1 \\ x_2 \end{pmatrix} = \begin{pmatrix} f_n \\ f_{n-1} \end{pmatrix}.
	\end{equation}
\end{lem}
\begin{proof}
	We first note d'Ocagne's Identity \cite[Ex. 357]{Gunderson} which asserts that
	\begin{equation} \label{Docc}
		f_{\alpha}f_{\beta + 1} - f_{\beta}f_{\alpha+1} = (-1)^\beta f_{\alpha-\beta}\quad\mbox{for all } \alpha,\beta\in \intg.
	\end{equation}
	Let $B_{k,\ell}$ be the $2\times 2$ matrix on the left hand side of \eqref{2x2Solve}.  Applying \eqref{Docc} with $\alpha = \ell - 1$ and $\beta = k -1$ we find that
	\begin{equation*}
		\det(B_{k,\ell} ) = f_{\ell-1} f_k - f_{k-1}f_\ell = (-1)^{k-1} f_{(\ell-1) - (k-1)} = (-1)^{k-1} f_{\ell-k}
	\end{equation*}
	which means that $B_{k,\ell}$ is invertible and \eqref{2x2Solve} has a unique solution.  Additionally, we know that
	\begin{equation*}
		B_{k,\ell}^{-1} = \frac{1}{(-1)^{k-1} f_{\ell - k}} \begin{pmatrix} f_{\ell-1} & - f_{\ell} \\ -f_{k-1} & f_k \end{pmatrix}
	\end{equation*}
	so that \eqref{2x2Solve} is equivalent to
	\begin{equation*}
		 \begin{pmatrix} x_1 \\ x_2 \end{pmatrix} 
		 	= \frac{1}{(-1)^{k-1} f_{\ell - k}}\begin{pmatrix} f_{\ell - 1}f_n - f_{\ell}f_{n-1} \\ f_k f_{n-1} - f_{k-1} f_n \end{pmatrix}.
	\end{equation*}
	By applying \eqref{Docc} in a similar manner as above, we obtain that
	\begin{equation*}
		 f_{\ell - 1}f_n - f_{\ell}f_{n-1} = (-1)^\ell f_{n-\ell}\quad\mbox{and}\quad f_k f_{n-1} - f_{k-1} f_n = (-1)^{k-1} f_{n-k}.
	\end{equation*}
	Therefore, \eqref{2x2Solve} is equivalent to 
	\begin{equation*}
		 \begin{pmatrix} x_1 \\ x_2 \end{pmatrix} = \frac{1}{(-1)^{k-1} f_{\ell - k}}\begin{pmatrix} (-1)^\ell f_{n-\ell} \\  (-1)^{k-1} f_{n-k} \end{pmatrix},
	\end{equation*}
	and we have established the lemma.
\end{proof}

For each $n\geq 1$ we define the $2\times n$ matrix
\begin{equation*}
	A_n =  \begin{pmatrix} f_1 & f_2 & \cdots & f_n \\ f_0 & f_1 & \cdots & f_{n-1} \end{pmatrix}
\end{equation*}
and note that $A_n$ defines a linear transformation from $\real^n$ to $\real^2$.  Clearly the rows of $A_n$ are linearly independent over $\real$ so that $\rank(A_n) = 2$ and $\dim(\ker(A_n)) = n-2$.  Let
\begin{equation*}
	\mathcal V_n(\real) = \left\{\xx\in \real^n:A_n(\xx) = \begin{pmatrix} f_n \\ f_{n-1} \end{pmatrix}\right\}
\end{equation*}
and note that $\mathcal V_n(\real)$ contains at least one element, namely $\ttt_n = (0,0,\ldots,0,1)^T\in \real^n$.
From our observation that $\rank(A_n) = 2$, $\mathcal V_n(\real)$ defines a set of points in $\real^n$ satisfying a pair of linear constraints.
In many situations, we will need to impose the additional constraint that $\xx$ has non-negative entries, and therefore, we also define
\begin{equation*}
	\mathcal V^+_n(\real) = \left\{ (x_1,x_2,\ldots,x_n)^T: x_i\geq 0\mbox{ for all } 1\leq i\leq n\right\}.
\end{equation*}
One easily checks that $\mathcal V_n^+(\real)$ is a compact subset of $\real^n$, and therefore, $\min\{g(\xx): \xx\in \mathcal V_n^+(\real)\}$ exists for all continuous functions $g:\real^n\to \real$.
We shall write
\begin{equation} \label{VnZ}
	\mathcal V_n^+(\intg) = \left\{ \xx\in \mathcal V^+_n(\real) : \xx\in \intg^n\right\}
\end{equation}
and note that $\mathcal V_n^+(\intg)$ is clearly a finite set.  The set $\mathcal V_n^+(\intg)$ appeared both in \cite{SamuelsExceptional} and in the introduction of the present paper where it was denoted
simply by $\mathcal V_n$, however, we find the more detailed notation \eqref{VnZ} to be more intuitive in the context our subsequent work.

We now describe an important strategy for minimizing a linear function of $\xx$ subject to the constraint that $\xx\in \mathcal V_n^+(\real)$.

\begin{lem} \label{ReductionToF}
	If $n$ is a positive integer then $\mathcal F_n^+ \subseteq \mathcal V_n^+(\real)$.  Moreover, if $L:\real^n \to \real$ is a linear map then the following conditions hold:
	\begin{enumerate}[label={(\roman*)}]
		\item \label{FVSwitch} $\min\{ L(\xx): \xx\in \mathcal V_n^+(\real)\} = \min\{ L(\xx): \xx\in \mathcal F_n^+\}$
		\item \label{IntegerMin} If there exists $\zz\in \intg^n$ which attains the minimum on the right hand side of \ref{FVSwitch}, then $\zz$ attains the minimum of the set $\{L(\xx): \xx\in \mathcal V_n^+(\intg)\}$.
		\item\label{UniqueMin} If the minimum on the right hand side of \ref{FVSwitch} is attained by a unique point $\zz\in \mathcal F_n^+$, then $\zz$ is the unique point attaining the minimum on
			the left hand side of \ref{FVSwitch}.
	\end{enumerate}
\end{lem}
\begin{proof}
	We first define the sets
	\begin{align*}
		\mathcal W_n  & =  \left\{(x_1,x_2,\ldots,x_n)^T\in \mathcal V_n(\real): x_i=0\mbox{ for all but at most two values of } i\right\} \\
		\mathcal W_n^+ & =  \left\{(x_1,x_2,\ldots,x_n)^T\in \mathcal W_n: x_i\geq 0\mbox{ for all } i\right\}
	\end{align*}
	and we claim that 
	\begin{equation} \label{WEqualities}
		\mathcal W_n^+ = \mathcal F_n^+.
	\end{equation}
	Assuming that $k,\ell\in \intg$ are such that $1\leq k < \ell\leq n$, we let $B_{k,\ell}$ be the matrix from the proof of Lemma \ref{2x2Lemma}.  From Lemma \ref{2x2Lemma}, we get that
	\begin{equation*}
		A_n\yy_n(k,\ell) = B_{k,\ell} \begin{pmatrix} (-1)^{\ell-k+1}f_{n-\ell}/f_{\ell - k} \\  f_{n-k}/f_{\ell - k} \end{pmatrix} =  \begin{pmatrix} f_n \\ f_{n-1} \end{pmatrix}
	\end{equation*}
	which shows that $\mathcal F_n\subseteq \mathcal W_n$.  On the other hand, if $\yy = (y_1,y_2,\ldots,y_n)^T \in \mathcal W_n$ then there exist $1\leq k < \ell \leq n$ such 
	that $y_i = 0$ for all $i\not\in \{k,\ell\}$.  This means that
	\begin{equation*}
		B_{k,\ell} \begin{pmatrix} y_k \\ y_\ell \end{pmatrix} = A_n\yy =  \begin{pmatrix} f_n \\ f_{n-1} \end{pmatrix}
	\end{equation*}
	and the uniqueness property of Lemma \ref{2x2Lemma} yields that $y_k = (-1)^{\ell-k+1}f_{n-\ell}/f_{\ell - k}$ and $y_\ell = f_{n-k}/f_{\ell - k}$.  We have now established that $\yy = \yy_n(k,\ell)$ and
	$\mathcal W_n \subseteq \mathcal F_n$, so we conclude that
	\begin{equation*}
		\mathcal W_n = \mathcal F_n.
	\end{equation*}
	Since $\mathcal W_n^+$ and $\mathcal F_n^+$ contain precisely those points of $\mathcal W_n$ and $\mathcal F_n$ whose entries are non-negative, respectively, we immediately obtain \eqref{WEqualities}. 
		
	Directly from the definition of $\mathcal W_n^+$ and \eqref{WEqualities}, we obtain that $\mathcal F_n^+\subseteq \mathcal V_n^+(\real)$.  Additionally, each vertex of $\mathcal V_n^+(\real)$ must belong to 
	$\mathcal W_n^+$, and as a result, both \ref{FVSwitch} and \ref{UniqueMin} follow from the Fundamental Theorem of Linear Programming (see \cite[\S 2.4]{LuenbergerYe} 
	or \cite[\S 3.5]{Vanderbei}, for example).  Finally, if $\zz\in \intg^n$ attains the minimum on the right hand side of \ref{FVSwitch} then $\zz\in \mathcal V_n^+(\intg)$.  We now obtain that 
	\begin{equation*}
		L(\zz) = \min\{L(\xx): \xx\in \mathcal F^+_n\} = \min\{ L(\xx): \xx\in \mathcal V_n^+(\real)\} \leq \min\{ L(\xx): \xx\in \mathcal V_n^+(\intg)\} \leq L(\zz)
	\end{equation*}
	completing the proof of \ref{IntegerMin}.
\end{proof}

While Lemma \ref{ReductionToF} describes a useful strategy for minimizing a linear function $\xx$ subject to certain constraints, it is not yet clear how this relates to our question on exceptional points.
Our next two lemmas establish the necessary connection.

\begin{lem} \label{UnifConvMeasureFunctions}
	Suppose that $n\in \nat$ and $\xx = (x_1,x_2,\ldots,x_n)$ is a vector of non-negative real numbers.  Further suppose that $\{\psi_k\}_{k=1}^\infty$ is a sequence of positive real numbers such
	that $\psi_k\to \phi$ as $k\to\infty$.  Further define $H_{\xx,k}:(0,\infty)\to \real$ by
	\begin{equation*}
		\quad H_{\xx,k}(t) = \sum_{j=1}^n x_i\max\{f_j,\psi_k f_{j-1}\}^t.
	\end{equation*}
	Then for every $T>0$, $H_{\xx,k}$ converges to $G_{\xx}$ uniformly on $(0,T]$ as $k\to \infty$.
\end{lem}
\begin{proof}
	If all entries of $\xx$ are equal to $0$, then clearly the result holds, so we assume that $x_i\ne 0$ for some $i$.  Let $\varepsilon > 0$.
	We must find $K\in\nat$ such that $|G_{\xx}(t) - H_{\xx,k}(t)| < \varepsilon$ for all $k\geq K$ and all $t\in (0,T]$.
	
	Assume without loss of generality that $n$ is even.  In this case, we have that
	\begin{equation} \label{FibInequalities}
		\frac{f_2}{f_1} < \frac{f_4}{f_3} < \cdots < \frac{f_n}{f_{n-1}} < \phi < \frac{f_{n-1}}{f_{n-2}} < \cdots < \frac{f_3}{f_2}
	\end{equation}
	Now choose $K_0\in \nat$ such that 
	\begin{equation*}
		\frac{f_n}{f_{n-1}} < \psi_k < \frac{f_{n-1}}{f_{n-2}}\quad\mbox{ for all } k\geq K_0.
	\end{equation*}
	Suppose that $k\geq K_0$.    For each odd positive integer $j$ we have
	\begin{equation*}
		\max\{f_j,\phi f_{j-1}\} = \max\{f_j,\psi_k f_{j-1}\} = f_j,
	\end{equation*}
	and for each even positive integer $j$ we get that
	\begin{equation*}
		\max\{f_j,\phi f_{j-1}\} = \phi f_{j-1}\quad \mbox{and}\quad \max\{f_j,\psi_k f_{j-1}\} = \psi_k f_{j-1}.	
	\end{equation*}
	Therefore,
	\begin{align*}
		|G_{\xx}(t) - H_{\xx,k}(t)| & = \left|  \sum_{j=1}^n x_i\left(\max\{f_j,\phi f_{j-1}\}^t - \max\{f_j,\psi_k f_{j-1}\}^t\right)\right| \\
			& = \left|\sum_{j=1}^n x_{2j} \left( (\phi f_{2j-1})^t - (\psi_k f_{2j-1})^t\right) \right| \\
			& = \left|\sum_{j=1}^n x_{2j}f_{2j-1}^t \left(\phi^t - \psi_k^t\right)\right|
	\end{align*}
	and we have established that
	\begin{equation} \label{GRewrite}
		|G_{\xx}(t) - H_{\xx,k}(t)| = |\phi^t - \psi_k^t|\cdot \left|\sum_{j=1}^n x_{2j}f_{2j-1}^t\right|.
	\end{equation}
	Clearly both expressions on the right hand side of \eqref{GRewrite} are increasing functions of $t$.  Now setting $M = \left|\sum_{j=1}^n x_{2j}f_{2j-1}^T\right|$ we obtain that
	\begin{equation*} \label{GInequality}
		|G_{\xx}(t) - H_{\xx,k}(t)| \leq M\cdot |\phi^T - \psi_k^T|.
	\end{equation*}
	Since $y\mapsto y^T$ is a continuous function on $[0,\infty)$ and $\psi_k$ converges to $\phi$, we know that $\psi_k^T$ converges to $\phi^T$.  Hence, we may choose $K \geq K_0$ such that 
	$|\phi^T - \psi_k^T| < \varepsilon/M$ for all $k\geq K$ and the result follows.
\end{proof}

While the value of $K$ described in the proof of Lemma \ref{UnifConvMeasureFunctions} is independent of $t$ (as it must be in order to establish uniform convergence), it does depend on both $T$ and $n$.
This means that any application of Lemma \ref{UnifConvMeasureFunctions} requires fixing those values in advance.  This fact will be evident in the proof of our next lemma.

\begin{lem} \label{MainV}
	Suppose that $n\geq 3$ is an integer and $t_1,t_2,\ldots,t_n,t_{n-1}$ are positive real numbers satisfying the following properties:
	\begin{enumerate}[label={(\roman*)}]
		\item\label{Decreasingts} $t_{i+1} < t_{i}$ for all $1\leq i\leq n-2$.
		\item\label{UniqueMinimum} For each $1\leq i\leq n-1$ there exists a unique point $\zz(i)\in \mathcal V_n^+(\intg)$ such that $G_{\zz(i)}(t_i) = \min\{G_{\xx}(t_i):\xx\in \mathcal V_n^+(\intg)\}$.
		\item\label{OneToOne} If $1\leq i\leq n-2$ we have that $\zz(i) \ne \zz(i+1)$.
	\end{enumerate}
	Then there exist infinitely many pairs of primes $(p,q)$ such that $p^{f_n}/q^{f_{n-1}}$ has at least $n-2$ exceptional points.
\end{lem}
\begin{proof}
	Suppose that $\{(p_k,q_k)\}_{k=1}^\infty$ is a sequence of pairs of primes such that 
	\begin{equation*}
		\lim_{k\to \infty}\frac{\log q_k}{\log p_k} = \phi.
	\end{equation*}
	We know that expressions of the form $\log q/\log p$, where $p$ and $q$ are prime, form a dense subset of $\real^+$, and therefore, such a sequence must always exist.  For simplicity, 
	we shall write $\psi_k = \log q_k/\log p_k$ so that $\psi_k$ satisfies the hypotheses of Lemma \ref{UnifConvMeasureFunctions}.
	
	Since $\zz(i)$ is the unique point in $\mathcal V_n^+(\intg)$ such that $G_{\zz(i)}(t_i) = \min\{G_{\xx}(t_i):\xx\in \mathcal V_n^+(\intg)\}$, we must have that $G_{\zz(i)}(t_i) < G_{\xx}(t_i)$ 
	for all $\xx\in \mathcal V_n^+(\intg)\setminus\{\zz(i)\}$.  Additionally, $\mathcal V_n^+(\intg)$ is a finite set so we may define
	\begin{equation*}
		\varepsilon = \min\left\{\frac{G_{\xx}(t_i) - G_{\zz(i)}(t_i)}{2} : 1\leq i\leq n-1\ \mbox{and}\ \xx\in \mathcal V_n^+(\intg)\setminus\{\zz(i)\} \right\} > 0.
	\end{equation*}
	Now let $H_{\xx,k}$ be as in the statement of Lemma \ref{UnifConvMeasureFunctions} and set and let $T = t_1+1$.  
	Applying Lemma \ref{UnifConvMeasureFunctions} and again using the fact that $\mathcal V_n^+(\intg)$ is a finite set, there exists $K\in \nat$ such that 
	\begin{equation} \label{UnifConv}
		\left| H_{\xx,k}(t) - G_{\xx}(t)\right| < \varepsilon
	\end{equation}
	for all $\xx\in \mathcal V_n^+(\intg)$, all $t\in (0,T]$, and all $k\geq K$.  
	
	We claim that
	\begin{equation} \label{HClaim}
		H_{\zz(i),k}(t_i) < H_{\xx,k}(t_i)
	\end{equation}
	for all $1\leq i\leq n-1$, all $\xx\in \mathcal V_n^+(\intg)\setminus\{\zz(i)\}$, and all $k\geq K$.  To see this, it follows from \eqref{UnifConv} that
	\begin{equation*}
		\left| H_{\zz(i),k}(t_i) - G_{\zz(i)}(t_i)\right| < \varepsilon,
	\end{equation*}
	which yields
	\begin{align*}
		H_{\zz(i),k}(t_i) & < G_{\zz(i)}(t_i) + \varepsilon \\
			& \leq  G_{\zz(i)}(t_i) + \frac{G_{\xx}(t_i) - G_{\zz(i)}(t_i)}{2} \\ 
			& = \frac{G_{\xx}(t_i) + G_{\zz(i)}(t_i)}{2}
	\end{align*}
	Similarly, \eqref{UnifConv} also implies that $\left| H_{\xx,k}(t_i) - G_{\xx}(t_i)\right| < \varepsilon$, from which we conclude that
	\begin{align*}
		H_{\xx,k}(t_i) & > G_{\xx}(t_i) - \varepsilon \\
			& \geq G_{\xx}(t_i) - \frac{G_{\xx}(t_i) - G_{\xx_n(i)}(t_i)}{2} \\
			& = \frac{G_{\xx}(t_i) + G_{\zz(i)}(t_i)}{2}. 
	\end{align*}
	We have now established \eqref{HClaim}.
	
	Let $\alpha_k = p_k^{f_n}/q_k^{f_{n-1}}$.  According to the main result of \cite{SamuelsContinued} (more succinctly summarized in \cite[Theorem 1.4]{SamuelsExceptional}), there exists $L>0$ such that
	\begin{equation} \label{MeasureDef}
		m_t(\alpha_k)^t = \min\left\{ \sum_{j=1}^n x_j m\left(\frac{p_k^{f_j}}{q_k^{f_{j-1}}}\right)^t: (x_1,x_2,\ldots,x_n)^T\in \mathcal V_n^+(\intg)\right\}
	\end{equation}
	for all $k\geq L$ and all $t>0$.  Now assuming that $k\geq \max\{K,L\}$, we obtain that
	\begin{align*}
		\sum_{j=1}^n x_j m\left(\frac{p_k^{f_j}}{q_k^{f_{j-1}}}\right)^t & = \sum_{j=1}^n x_j\max\left\{ \log p_k^{f_j}, \log q_k^{f_{j-1}}\right\}^t \\
			& = \sum_{j=1}^n x_j\max\left\{ f_j\log p_k, f_{j-1} \log q_k\right\}^t \\
			& = (\log p_k)^t  \sum_{j=1}^n x_i\max\left\{ f_j, \psi_k f_{j-1} \right\}^t \\
			& = (\log p_k)^t H_{\xx,k}(t),
	\end{align*}
	and \eqref{MeasureDef} leads to 
	\begin{equation} \label{MeasureH}
		\frac{m_t(\alpha)}{\log p_k} = \left(\min\left\{ H_{\xx,k}(t): \xx\in \mathcal V_n^+(\intg)\right\}\right)^{1/t}.
	\end{equation}
	We now fix $i\in \intg$ with $1\leq i\leq n-2$.  Because of \eqref{HClaim} and condition \ref{OneToOne}, the minima of the sets $\left\{ H_{\xx,k}(t_i): \xx\in \mathcal V_n^+(\intg)\right\}$ and 
	$\left\{ H_{\xx,k}(t_{i+1}): \xx\in \mathcal V_n^+(\intg)\right\}$  must be attained by two distinct points in $\mathcal V_n^+(\intg)$.  In other words, the minimum on the right hand side of \eqref{MeasureH} must be attained by 
	distinct points when $t = t_i$ and $t=t_{i+1}$.  Therefore, we may apply \cite[Theorem 1.4]{SamuelsParametrized} to conclude that $(t_{i+1},t_{i})$ must contain at least one exceptional point for all $1\leq i\leq n-2$.  
	Hence, $\alpha_k$ must have at least $n-2$ exceptional points and the result follows.
\end{proof}

The proof of Theorem \ref{Main} is now a simple application of Lemmas \ref{ReductionToF} and \ref{MainV}.

\begin{proof}[Proof of Theorem \ref{Main}]
	We need only verify that the assumptions of Lemma \ref{MainV} hold.  Certainly \ref{Decreasingts} and \ref{OneToOne} follow directly from the assumptions of Theorem \ref{Main}, so it remains only to
	establish \ref{UniqueMinimum}, and for this purpose, we fix $i\in \intg$ such that $1\leq i\leq n-1$.  
	
	By applying \ref{UniquenessMain}, there exists a unique point $\zz\in \mathcal F_n^+$ such that $G_{\zz}(t_i) = \min\{ G_{\xx}(t_i):\xx\in \mathcal F^+_n\}$.
	Since $\xx\mapsto G_{\xx}(t)$ is a linear map, it follows from Lemma \ref{ReductionToF}\ref{FVSwitch} that 
	\begin{equation} \label{GVR}
		G_{\zz}(t_i) = \min\{ G_{\xx}(t_i):\xx\in \mathcal V_n^+(\real)\}.
	\end{equation}
	Then using Lemma \ref{ReductionToF}\ref{IntegerMin} and \ref{IntegersMain}, we obtain that
	\begin{equation} \label{GV}
		G_{\zz}(t_i) = \min\{ G_{\xx}(t_i):\xx\in \mathcal V_n^+(\intg)\}.
	\end{equation}
	If there is another point $\yy\in \mathcal V_n^+(\intg)$ such that $G_{\yy}(t_i) = G_{\zz}(t_i)$, then \eqref{GVR} implies that $$G_{\yy}(t_i) = \min\{ G_{\xx}(t_i):\xx\in \mathcal V_n^+(\real)\}$$
	contradicting Lemma \ref{ReductionToF}\ref{UniqueMin}.  It now follows that $\zz$ is the unique point in $\mathcal V_n^+(\intg)$ satisfying \eqref{GV}.  This verifies Lemma \ref{MainV}\ref{UniqueMinimum},
	and completes the proof of Theorem \ref{Main}.
\end{proof}

\section{Proof of Theorem \ref{MainExceptional}} \label{Application}

Our proof of Theorem \ref{MainExceptional} is purely a computation performed in MATLAB designed to verify the assumptions of Theorem \ref{Main}.
Before proceeding with that computation, we must select values of $t_i$.  Recall that there exists a unique positive real number $s_i$ such that
\begin{equation} \label{si2}
	\max\{f_{i+2},\phi f_{i+1}\}^{s_i} = \max\{f_{i+1},\phi f_{i}\}^{s_i} + \max\{f_{i},\phi f_{i-1}\}^{s_i},
\end{equation}
and moreover, $\{s_i\}_{i=1}^\infty$ is a strictly decreasing sequence with $s_i>1$ for all $i$.  Based on evidence presented in \cite{SamuelsExceptional}, we pose the following conjecture.

\begin{conj} \label{MainConjecture}
	Assume that $n\geq 3$ is an integer and $\{t_i\}_{i=1}^\infty$ is a sequence satisfying
	\begin{equation} \label{ti2}
		t_1 > s_1 > t_2 > s_2 > t_3 > s_3 > \cdots.
	\end{equation}
	Then for every $1\leq i\leq n-1$, $\yy_n(i,i+1)$ is the unique element of $\mathcal F_n^+$ such that
	\begin{equation} \label{YMin}
		G_{\yy_n(i,i+1)}(t_i) = \min\left\{G_{\xx}(t_i):\xx\in \mathcal F_n^+\right\}.
	\end{equation}
\end{conj}

If Conjecture \ref{MainConjecture} holds, then clearly the assumptions of Theorem \ref{Main} are satisfied regardless of the value of $n$.  Indeed, we certainly have that $\yy_n(i,i+1) \in \intg^n$ and 
$\yy_n(i,i+1) \ne \yy_n(i+1,i+2)$ for all relevant values of $i$.  In this scenario, we would obtain a proof of Conjecture \ref{ExceptionalCounting}.  
These observations suggest that any sequence $\{t_i\}_{i=1}^\infty$ which satisfies \eqref{ti2} will also satisfy the assumptions of Theorem \ref{Main}.
Using the values of $t_i$ in Figure \ref{SValues}, we may execute the following algorithm for each $1\leq i\leq n-1$:
\begin{enumerate}
	\item\label{Calcs} Calculate $G_{\xx}(t_i)$ for each $\xx\in \mathcal F_n^+$
	\item\label{Minimize} Calculate the value of $\xx$ which attains the minimum  of the values in \eqref{Calcs}
	\item Verify that the value of $\xx$ from \eqref{Minimize} is unique
	\item Verify that the value of $\xx$ from \eqref{Minimize} equals $\yy_n(i,i+1)$
\end{enumerate}
Because of \eqref{FTCardinality}, this algorithm requires only polynomial time and is completed easily for $n=39$.  We performed this work in MATLAB verifying that \eqref{YMin} holds for all $1\leq i\leq 38$.
Since $\yy_n(i,i+1) \in \intg^n$ and $\yy_n(i,i+1) \ne \yy_n(i+1,i+2)$ for all relevant values of $i$, this verifies that the assumptions of Theorem \ref{Main} hold when $n=39$.  These observations
complete the proof of Theorem \ref{MainExceptional}.

When we reach $n=40$ our algorithm fails because of a machine precision error in MATLAB.  Indeed, there exist two distinct points 
$\xx,\yy\in \mathcal F_{40}^+$ for which MATLAB is unable to distinguish between $G_{\xx}(t_{38})$ and $G_{\yy}(t_{38})$, and therefore, we are unable to establish the uniqueness property of Theorem \ref{Main}. 
Note that this type of failure is distinct from that which occurs in \cite{SamuelsExceptional}.  As we have noted, that work used an exponential time algorithm which failed when $n=14$ because its running time was too long.
Our present methods use the polynomial time algorithm described above, and hence, running time is of no concern for these small values of $n$.

While unfortunate, our machine precision error is completely expected.  Indeed, the values of $s_i$ are defined so that $\{G_{\xx}(s_i):\xx\in \mathcal F_n^+\}$ {\bf does not} have a unique minimum for any $1\leq i\leq n-2$.
Moreover, we must always have that $s_{n-2} < t_{n-2} < s_{n-3}$, meaning that $t_{n-2}$ is at most $(s_{n-3} - s_{n-2})/2$ units from a point $s$ where $\{G_{\xx}(s):\xx\in \mathcal F_n^+\}$
has no unique minimum.  Since
\begin{equation*}
	\frac{s_{37} - s_{38}}{2} = 1.36011\times 10^{-16}
\end{equation*}
and MATLAB's machine epsilon is $2.2204\times 10^{-16}$, our methods will necessarily fail when $n=40$.  This value of machine epsilon results from using the IEEE standard double-precision floating-point format. 
Although it may be possible to extend this using software packages with arbitrarily high precision arithmetic, these methods can use an approximation to the field of real numbers, and may yield 
answers which differ from the true answer. 

\begin{figure}
	$\begin{array}{c | c | c}
		i & s_i -1 & t_i -1\\ \hline
		1 & 7.146171809\times 10^{-1} & 8\times 10^{-1} \\
		2  & 1.940660944\times 10^{-1} &  2\times 10^{-1} \\
		3 & 7.478824327\times10^{-2} & 8\times 10^{-2} \\ 
		4 & 2.743232714\times 10^{-2}& 3\times 10^{-2} \\
		5 & 1.050775991\times 10^{-2}& 2\times 10^{-2} \\
		6 & 3.990473545\times 10^{-3} & 4\times 10^{-3} \\
		7 & 1.524905897\times 10^{-3}& 2\times 10^{-3} \\
		8 & 5.819727416\times 10^{-4} & 6\times 10^{-4} \\
		9 & 2.223084584\times 10^{-4} &  3\times 10^{-4}\\
		10 & 8.490386498\times 10^{-5} & 9\times 10^{-5} \\
		11 & 3.243070303\times 10^{-5} & 4\times 10^{-5} \\
		12 & 1.238720471\times 10^{-5} & 2\times 10^{-5} \\
		13 & 4.731497826\times 10^{-6} & 5\times 10^{-6} \\
		14 & 1.807266635\times 10^{-6} & 2\times 10^{-6} \\
		15 & 6.903145697\times 10^{-7} & 7\times 10^{-7}  \\
		16 & 2.636766023\times 10^{-7} & 3\times 10^{-7} \\
		17 & 1.007155030\times 10^{-7} & 2\times 10^{-7} \\
		18 & 3.846989684\times 10^{-8} & 4\times 10^{-8}  \\
		19 & 1.469419311\times 10^{-8} & 2\times 10^{-8} \\
		20 & 5.612682286\times 10^{-9} & 6\times 10^{-9} \\
		21 & 2.143853866\times 10^{-9} & 3\times 10^{-9} \\
		22 & 8.188793092\times 10^{-10} & 9\times 10^{-10} \\
		23 & 3.127840634\times 10^{-10} & 4\times 10^{-10}  \\
		24 & 1.194728810\times 10^{-10} & 2\times 10^{-10} \\
		25 & 4.563457984\times 10^{-11}&  5\times 10^{-11} \\
		26 & 1.743085843\times 10^{-11} & 2\times 10^{-11} \\
		27 & 6.657995469\times 10^{-12} & 7\times 10^{-12} \\
		28 & 2.543127972\times 10^{-12}& 3\times 10^{-12} \\
		29 & 9.713884477\times 10^{-13} & 1\times 10^{-12} \\
		30 & 3.710373707\times 10^{-13} & 4\times 10^{-13} \\
		31 & 1.417236645\times 10^{-13}& 2\times 10^{-13} \\
		32 & 5.413362284\times 10^{-14} & 6\times 10^{-14} \\
		33 & 2.067720399\times 10^{-14} & 3\times 10^{-14} \\
		34 & 7.897989132\times 10^{-15} & 9\times 10^{-15} \\
		35 & 3.016763405\times 10^{-15} & 5\times 10^{-15} \\
		36 & 1.152301085\times 10^{-15} & 2\times 10^{-15} \\
		37 & 4.401398491\times 10^{-16} & 8\times 10^{-16}\\
		38 & 1.681184625\times 10^{-16} & 2\times 10^{-16} \\
	\end{array}$
	\caption{Approximate values of $s_i$ as defined by \eqref{si2} with our corresponding choices of $t_i$}
	\label{SValues}
\end{figure}


\section{Supplementary Results} \label{Supplementary}

We noted in the introduction that the cardinality of $\mathcal V_n^+(\intg)$ (denoted by $\mathcal V_n$ at that time) grows at least exponentially as a function of $n$.  We now provide the reader with a formal proof of this assertion.

\begin{thm} \label{ExpGrowth}
	If $n\geq 3$ then $\#\mathcal V_n^+(\intg) \geq f_n$.
\end{thm}
\begin{proof}
	We first define vectors $\bb_n(1),\bb_n(2),\ldots,\bb_n(n-2)\in \real^n$ by
	\begin{equation} \label{BVectors}
		\bb_n(1) = \begin{pmatrix} -1 \\ -1 \\ 1 \\ 0 \\ 0 \\ \vdots \\ 0 \end{pmatrix},\ \bb_n(2) = \begin{pmatrix} 0 \\ -1 \\ -1 \\ 1 \\ 0 \\ \vdots \\ 0 \end{pmatrix},\ \cdots,\ 
			\bb_n(n-2) = \begin{pmatrix} 0 \\ \vdots \\ 0 \\ 0 \\ -1 \\ -1 \\ 1 \end{pmatrix}.
	\end{equation}
	By the recurrence relation in the Fibonacci sequence, we conclude that $\bb_i\in \ker(A_n)$ for all $1\leq i\leq n-2$.  For each such $i$, we also define
	\begin{equation*}
		\Lambda_n(i) = \left\{ \yy_n(i,i+1) +k\bb_n(i):0\leq k < f_{n-i-1}\right\}
	\end{equation*}
	and note the following facts regarding $\Lambda_n(i)$:
	\begin{enumerate}[label={(\roman*)}]
		\item $\Lambda_n(i) \subseteq \mathcal V_n^+(\intg)$
		\item $\Lambda_n(i) \cap \Lambda_n(j) = \emptyset$ for all $i\ne j$
		\item $\#\Lambda_n(i) = f_{n-i-1}$
		\item $\yy_n(n-1,n)$ does not belong to $\Lambda_n(i)$ for any $1\leq i\leq n-2$.
	\end{enumerate}
	We now conclude that
	\begin{equation*}
		\#\mathcal V_n^+(\intg) \geq 1+  \sum_{i=1}^{n-2} f_{n-i-1} = 1+ \sum_{i=1}^{n-2} f_i = 1 + f_n - 1 = f_n.
	\end{equation*}
\end{proof}
Since it is well-known that $f_n/\phi^n \to 1/\sqrt 5$ as $n\to \infty$, it follows from Theorem \ref{ExpGrowth} that $\#\mathcal V_n^+(\intg)$ exhibits at least exponential growth.

The primary advantage of Theorem \ref{Main} over the relevant previous work \cite{SamuelsExceptional} is that it improves a search algorithm for verifying the existence of rational numbers with at
least $n-2$ exceptional points.  We now verify that we do indeed obtain the advertised improvement.  Specifically, we now show how to compute the cardinality of $\mathcal F_n^+ \times \mathcal T_n$.

\begin{thm} \label{FCard}
	If $n\geq 3$ is an integer then we have the following formulas:
	\begin{align*}
		\#\mathcal F_n & = \frac{n^2 - 3n+4}{2} \\
		\#\mathcal F^+_n & = \begin{cases} (n^2 - 2n + 5)/4 & \mbox{if } n\mbox{ is odd} \\ (n^2 - 2n + 4)/4 & \mbox{if } n\mbox{ is even}. \end{cases} \\
		\#(\mathcal F^+_n\times\mathcal T_n) & = \begin{cases} (n^3 -3n^2 + 7n - 5)/4 & \mbox{if } n\mbox{ is odd} \\ (n^3 - 3n^2 + 6n - 4)/4 & \mbox{if } n\mbox{ is even}. \end{cases}
	\end{align*}
\end{thm}
\begin{proof}
	If $\ell = n$ then $\yy_n(k,\ell) = \ttt_n$ for all $1\leq k < \ell$, and therefore
	\begin{equation} \label{FSplit}
		\mathcal F_n =  \left\{\yy_n(k,\ell): 1\leq k < \ell \leq n-1\right\} \cup \{\ttt_n\}.
	\end{equation}
	Moreover, if $\yy_n(k,\ell)$ is such that $1\leq k < \ell \leq n-1$ and $\yy_n(k,\ell) = (y_1,y_2,\ldots,y_n)^T$ then $y_i\ne 0$ if and only if $i\in \{k,\ell\}$.
	Consequently, 
	\begin{equation*}
		\#\left\{\yy_n(k,\ell): 1\leq k < \ell \leq n-1\right\} = \#\{(k,\ell): 1\leq k < \ell \leq n-1\} = {n-1\choose 2}
	\end{equation*}
	and it follows from \eqref{FSplit} that
	\begin{equation*}
		\#\mathcal F_n = {n-1\choose 2} + 1 = \frac{n^2 - 3n+4}{2}
	\end{equation*}
	establishing the first formula.  
	
	To verify second formula, we once again note that
	\begin{equation} \label{F+Split}
		\#\mathcal F^+_n =  \left\{\yy_n(k,\ell): 1\leq k < \ell \leq n-1,\ \ell - k\ \mathrm{is\ odd}\right\} \cup \{\ttt_n\}
	\end{equation}
	and
	\begin{equation} \label{Translate}
		\#\left\{\yy_n(k,\ell): 1\leq k < \ell \leq n-1,\ \ell - k\ \mathrm{is\ odd}\right\} = \#\left\{(k,\ell):1\leq k < \ell \leq n-1,\ \ell - k\ \mathrm{is\ odd}\right\}.
	\end{equation}	
	Clearly $\ell - k$ is even if and only if both $k$ and $\ell$ are even or both $k$ and $\ell$ are odd.  Therefore, the right hand side of \eqref{Translate} is equal to
	\begin{equation*}
		{n-1\choose 2} - \#\left\{(k,\ell):1\leq k < \ell \leq n-1,\ k,\ell\ \mathrm{are\ even}\right\} - \#\left\{(k,\ell):1\leq k < \ell \leq n-1,\ k,\ell\ \mathrm{are\ odd}\right\}.
	\end{equation*}
	If $n$ is odd then there are $(n-1)/2$ even integers and $(n-1)/2$ odd integers in $[1,n-1]$.  Similarly, if $n$ is even then the are $n/2$ odd integers and $(n-2)/2$ even integers in $[1,n-1]$.
	For odd $n$, we conclude that
	\begin{equation*}
		\#\left\{(k,\ell):1\leq k < \ell \leq n-1,\ \ell - k\ \mathrm{is\ odd}\right\} = {n-1 \choose 2} - 2\cdot {(n-1)/2\choose 2},
	\end{equation*}
	and for even $n$, we get that
	\begin{equation*}
		\#\left\{(k,\ell):1\leq k < \ell \leq n-1,\ \ell - k\ \mathrm{is\ odd}\right\} = {n-1 \choose 2} - {(n-2)/2\choose 2} - {n/2 \choose 2}.
	\end{equation*}
	These expressions simplify to $(n^2 - 2n + 1)/4$ and $(n^2 - 2n)/4$, respectively, and the second formula follows from \eqref{F+Split} and \eqref{Translate}. 
	
	Finally, since $t_3,t_4,\ldots,t_{n+1}$ is a strictly decreasing sequence, we know that $\#\mathcal T_n = n-1$.  Hence, 
	$\#(\mathcal F^+_n\times\mathcal T_n)  = (\#\mathcal F_n^+) \cdot (\#\mathcal T_n)$ and the final formula follows from a basic simplification.
\end{proof}


\end{document}